\documentclass[11pt]{article}
\usepackage{array, amsmath, amssymb}

\begin{document}
\begin{center}
{\Large\bf Duality and upper bounds in optimal stochastic control \\
\vspace{0.5cm}
governed by partial differential equations}\\
\vspace{0.5cm}
{\bf Shinji Tanimoto}\\ 
\vspace{0.5cm}
Department of Mathematics,
University of Kochi,\\
Kochi 780-8515, Japan\footnote{Former affiliation}.  \\
\end{center}
\begin{abstract}
A dual control problem is presented for the optimal stochastic control of a system governed by partial differential 
equations. Relationships between the optimal values of the original and the dual problems are investigated and two
duality theorems are proved. The dual problem serves to provide upper bounds for the optimal and maximum value
of the original one or even to give the optimal value. 
\end{abstract}
\vspace{0.7cm}
{\large \bf 1. Introduction} \\  
\\
\indent
The original problem (or primal problem) considered is the optimal control of a system governed by a
stochastic heat equation that is described in [4], which is a maximization problem. In this paper,
to the problem we associate another, called its {\it dual problem}, which is in turn a minimization problem. 
We prove two types of duality theorem. \\
\indent
First we show that solutions of the dual problem provide upper bounds for the maximum of the primal problem. We call 
this assertion a {\it weak duality theorem}. Next,
under some conditions related to the maximum principle of control theory, 
the maximum can be attained by solving the dual problem. Such a property is called a {\it strong duality theorem}. \\
\indent
Let $T > 0$ and $V$ be a bounded and open domain in $\mathbb{R}^n$ with $C^1$ boundary $\partial V = \Gamma$. 
On $[0, T] \times V$ we consider the following stochastically controlled system. 
The one-dimensional Brownian motion $B(t) = B(t, \omega)$ is defined on a filtered probability space 
$(\Omega, \mathcal{F}, \{\mathcal{F}_t \}_{t \ge 0}, P)$. The state of the system is denoted by 
$X(t, x) \in \mathbb{R}$, which is controlled by $u(t, x) \in \mathbb{R}$ for $t \in [0, T]$ and at 
$x \in \bar{V} = V \cup \Gamma$. The control process $u(t, x) = u(t, x, \omega)$ satisfies 
$u(t, x) \in \mathcal{U}$, where $\mathcal{U}$ is a bounded set of $\mathbb{R}^k$,
and it is $\mathcal{F}_t$-measurable for all $(t, x) \in (0, T) \times V$.  
The state $X(t, x)$ is  
described by a stochastic heat equation of the form
\begin{eqnarray}
dX(t, x) & = & \big(\hat{A}X(t, x) + C(t, x, u(t, x))\big)dt + \sigma(t, x) dB(t), \\
X(0, x) & =& \xi(x)~ ~{\rm for}~x \in \bar{V}, \\
X(t, x) & = & \eta(t, x)~~{\rm for}~(t, x) \in (0, T) \times \Gamma. 
\end{eqnarray}
The boundary value functions $\xi$ on $\bar{V}$, and $\eta$ on $(0, T) \times \Gamma$ are $C^1$ real-valued and deterministic.
$\hat{A}$ is a second order partial differential operator acting on smooth functions of $x$:
\begin{eqnarray*}
\hat{A}\phi (x) = \sum_{i, j = 1}^{n} a_{ij}(x) \frac{{\partial}^2 \phi(x)}{\partial x_i \partial x_j} + 
\sum_{i = 1}^{n} b_i(x) \frac{\partial \phi(x)}{\partial x_i},
\end{eqnarray*}
where $(a_{ij}(x))$ is a symmetric nonnegative definite $n \times n$ matrix with entries 
$a_{ij}(x) \in C^2(V) \cap C(\bar{V})$ and $b_i(x) \in C^2(V) \cap C(\bar{V})$ for $1 \le i \le n$. 
A control process $u(t, x)$ is called {\it admissible} 
if the corresponding solution $X_u(t, x)$ of  Eqs.(1)-(3) is unique and belongs to $L^2(\Lambda \times P)$,
where $\Lambda$ is the Lebesgue measure on $[0, T] \times \bar{V}$. 
The set of all admissible controls is denoted by $\mathcal{A}$;
\begin{eqnarray*}
\mathcal{A} = \{ u(t, x) ~|~ u(t, x) \in \mathcal{U} ~{\rm is}~ \mathcal{F}_t {\rm -measurable ~for~ all} ~(t, x) \}.
\end{eqnarray*}
The $C^1$ functions $C$ and $\sigma$ in (1) are, respectively,  $C: [0, T] \times V \times \mathcal{U} \to \mathbb{R}$
and $\sigma: [0, T] \times V \to \mathbb{R}$.
\indent
The expected performance (or payoff) is given by, for each $u \in \mathcal{A}$, 
\begin{eqnarray}
{\mathcal J}(u) = \mathbb{E} \Big[ \int_0^T \Big( \int_V \big(F(t, x, X(t, x)) + G(t, x, u(t,x)) \big)dx \Big)dt \Big],
\end{eqnarray}
where $X(t, x) = X_u(t, x)$. Throughout this paper we impose the following; \\
\\
(Assumption) $F: [0, T] \times V \times \mathbb{R} \to \mathbb{R}$ is a $C^1$ function that is 
{\it concave} with respect to $X$.  \\
\\
$G: [0, T] \times V \times \mathcal{U} \to \mathbb{R}$ is a bounded continuous function, and
$\mathbb{E}$ denotes the expectation with respect to the probability measure $P$. 
The aim of the primal problem is to find a maximizing control $u^{\ast} \in \mathcal{A}$ and ${\mathcal J}^{\ast} \in \mathbb{R}$ 
such that
\begin{eqnarray*}
{\mathcal J}^{\ast} = \sup_{u \in \mathcal{A}} {\mathcal J}(u) = {\mathcal J}(u^{\ast}).
\end{eqnarray*}
Thus the primal problem is formulated as
\begin{eqnarray}
\sup_{u \in \mathcal{A}} {\mathcal J}(u).
\end{eqnarray}
\indent
In the next section a dual problem to (5) is proposed. 
Similar dual control problems were constructed for max-min control problems in [5], for
non-well-posed distributed systems in [6] and for optimal stochastic control in [7]. When a primal problem is a 
minimization problem, its dual problem serves to provide lower bounds for the minimum value of the primal one. 
Here the primal problem is a maximization problem, its dual problem provides upper bounds for the maximum value. 
Under some conditions related to the maximum principle of control theory it is also able to attain the maximum. \\
\\
\\
{\large \bf 2. Dual Problem} \\  
\\
\indent 
The adjoint of the differential operator $\hat{A}$ is defined by
\begin{eqnarray*}
\hat{A}^{\ast}\phi (x) = \sum_{i, j = 1}^{n} \frac{{\partial}^2 (a_{ij}(x)\phi(x))}{\partial x_i \partial x_j} -
\sum_{i = 1}^{n} \frac{\partial (b_i(x)\phi(x))}{\partial x_i}.
\end{eqnarray*}
In order to present the dual problem, for each real number $p$, we define a function
\begin{eqnarray}
H(t, x, p) = \sup_{u \in \mathcal U} \,\big(G(t, x, u) + p\,C(t, x, u)\big).
\end{eqnarray}
Note that the control variable $u$ of the primal problem disappears at this stage.  \\
\indent
The dual control problem is the system with performance functional that is to be minimized:
\begin{eqnarray}
\mathcal{L}(X, p) = \mathbb{E} \Big[ \int_V \int_0^T \Big( H(t, x, p(t, x)) + F(t, x, X(t, x)) - 
X(t, x) \frac{\partial F(t, x, X(t, x))}{\partial X}\Big) dtdx  \nonumber\\
- \int_V \int_0^T \Big(X(t, x)\hat{A}^{\ast}p(t, x) - p(t, x)\hat{A}X(t, x) \Big)dtdx  \\
+ \int_V \Big(p(0, x)\xi(x) + \int_0^T p(t, x)\sigma(t, x)dB(t)\Big) dx \Big],    \nonumber
\end{eqnarray}
over all variables $X$ and $p$ that satisfy:
\begin{eqnarray}
- dp(t, x) /dt = \hat{A}^{\ast}p(t, x) + \partial F(t, x, X(t, x)) /\partial X ; ~~~~~~~~~~\\
p(T, x) = 0 ~~{\rm for}~x \in \bar{V},  ~~~p(t, x) = 0 ~~{\rm for}~(t, x) \in (0, T) \times \Gamma;  ~~~  \\
X(0, x) = \xi(x) ~~{\rm for}~x \in \bar{V}, ~~~X(t, x) = \eta(t, x)~~{\rm for}~(t, x) \in (0, T) \times \Gamma. 
\end{eqnarray}
The variable $X(t, x)$ plays a role of control process of the dual problem 
that is a continuous process belonging to $L^2(\Lambda \times P)$. 
As indicated by the strong duality theorem (Section 4), $X(t, x)$ may be a solution of Eqs.(1)--(3),
which indeed becomes a continuous process. Or it can be even a deterministic and continuous variable.
Hence the dual problem is more manageable than the primal one. 
The variable $p(t, x)$, in turn, represents the state of the dual problem. 
We denote by $\mathcal B$ the set of all pairs $(X, p)$ that satisfy Eqs.(8)--(10).
So the dual problem is formulated as
\begin{eqnarray}
\inf_{(X,~p) \in \mathcal B} \mathcal{L}(X, p).
\end{eqnarray}
\\
\\
{\large \bf 3. Weak Duality Theorem} \\  
\\
\indent 
In this section we show that solutions of the dual problem provide upper bounds for the maximum of problem (5). 
We call this property a {\it weak duality theorem}.  \\
\\
{\bf Theorem 1.} {\it Under the concavity of the function $F$ it follows that}
\begin{eqnarray*}
   \sup_{u \in \mathcal A} {\mathcal J}(u) \le \inf_{(X,~p) \in {\mathcal B}}{\mathcal L}(X, p).
\end{eqnarray*}
\\
{\it Proof.} 
Let $u \in \mathcal{A}$ be an admissible control and let us fix it for the moment. Let $\bar{X}$ be the solution of Eqs.(1)-(3)
for $u$ and put (see Eq.(4));
\begin{eqnarray*}
\mathcal J(u) = \mathbb{E} \Big[ \int_0^T \Big( \int_V \big(F(t, x, \bar{X}(t, x)) + G(t, x, u(t,x)) \big)dx \Big)dt \Big].
\end{eqnarray*}
On the other hand, for the same $u$ we consider the following expectation,
 using an arbitrary $(X^{\circ}, p^{\circ}) \in \mathcal{B}$:
\begin{eqnarray}
L(u; X^{\circ}, p^{\circ}) = 
\mathbb{E} \Big[ \int_V \int_0^T \Big(F(t, x, X^{\circ}(t, x)) + G(t, x, u(t, x)) + p^{\circ}(t, x)C(t, x, u(t, x)) \\
- X^{\circ}(t, x) \frac{\partial F(t, x, X^{\circ}(t, x))}{\partial X}\Big) dtdx  
- \int_V \int_0^T \Big( X^{\circ}(t, x)\hat{A}^{\ast}p^{\circ}(t, x) - p^{\circ}(t, x)\hat{A}X^{\circ}(t, x) \Big)dtdx  \nonumber\\
+ \int_V \Big(p^{\circ}(0, x)\xi(x) + \int_0^T p^{\circ}(t, x)\sigma(t, x)dB(t)\Big) dx \Big],  \nonumber
\end{eqnarray}
where $(X^{\circ}, p^{\circ})$ is a solution of Eqs.(8)-(10):
\begin{eqnarray*}
- dp^{\circ}(t, x)/dt &=& \hat{A}^{\ast}p^{\circ}(t, x) + \partial F(t, x, X^{\circ}(t, x)) /\partial X; \\
p^{\circ}(T, x) &=& 0 ~~{\rm for}~x \in \bar{V}, ~~    
p^{\circ}(t, x) = 0 ~~{\rm for}~(t, x) \in (0, T) \times \Gamma;  \\
X^{\circ}(0, x) &=& \xi(x) ~~{\rm for}~x \in \bar{V}, ~~
X^{\circ}(t, x) = \eta(t, x)~~{\rm for}~(t, x) \in (0, T) \times \Gamma. 
\end{eqnarray*}
Making use of these fixed $u \in \mathcal A$ and $(X^{\circ}, p^{\circ}) \in \mathcal{B}$, the difference between 
${\mathcal J}(u)$ and $L(u; X^{\circ}, p^{\circ})$ is
\begin{eqnarray*}
~~~~~{\mathcal J}(u) - L(u; X^{\circ}, p^{\circ}) =  
\mathbb{E} \Big[ \int_V \int_0^T \Big(F(t, x, \bar{X}(t, x)) - F(t, x, X^{\circ}(t, x)) - p^{\circ}(t, x)C(t, x, u(t, x)) \\
+  X^{\circ}(t, x) \frac{\partial F(t, x, X^{\circ}(t, x))}{\partial X}\Big) dtdx 
+ \int_V \int_0^T \Big( X^{\circ}(t, x)\hat{A}^{\ast}p^{\circ}(t, x) - p^{\circ}(t, x)\hat{A}X^{\circ}(t, x) \Big)dtdx \\
- \int_V \Big(p^{\circ}(0, x)\xi(x) + \int_0^T p^{\circ}(t, x)\sigma(t, x)dB(t)\Big) dx \Big].
\end{eqnarray*}
By the concavity of $F$ with respect to $X$ we have
\begin{eqnarray*}
F(t, x, \bar{X}(t, x)) - F(t, x, X^{\circ}(t, x)) \le \frac{\partial F(t, x, X^{\circ}(t, x))}{\partial X}(\bar{X}(t, x) - 
X^{\circ}(t, x)),
\end{eqnarray*}
from which we have the inequality 
\begin{eqnarray}
{\mathcal J}(u) - L(u; X^{\circ}, p^{\circ}) \le 
\mathbb{E} \Big[ \int_V \int_0^T \Big (\frac{\partial F(t, x, X^{\circ}(t, x))}{\partial X}\bar{X}(t, x)
- p^{\circ}(t, x)C(t, x, u(t, x)) \Big) dtdx  \nonumber \\
+  \int_V \int_0^T \Big(X^{\circ}(t, x)\hat{A}^{\ast}p^{\circ}(t, x) - p^{\circ}(t, x)\hat{A}X^{\circ}(t, x) \Big)dtdx  
~~~~~~~~~~~~~\\
- \int_V \Big(p^{\circ}(0, x)\xi(x) + \int_0^T p^{\circ}(t, x)\sigma(t, x)dB(t)\Big) dx \Big].  \nonumber
\end{eqnarray}
We show that the right-hand side of (13) is equal to zero. From Eq.(8) it follows that
\begin{eqnarray*} \partial F(t, x, X^{\circ}(t, x))/\partial X = - dp^{\circ}(t, x) / dt - \hat{A}^{\ast}p^{\circ}(t, x),
\end{eqnarray*}
and that
\begin{eqnarray}
{\mathcal J}(u) - L(u; X^{\circ}, p^{\circ}) \le -
\mathbb{E} \Big[ \int_V \int_0^T \Big(\frac{dp^{\circ}(t, x)}{dt}\bar{X}(t, x) + 
\bar{X}(t, x)\hat{A}^{\ast}p^{\circ}(t, x)  \nonumber ~~~~~~~~~~~~~~~~~~~~~~~ \\
~~~~~~~~ + p^{\circ}(t, x)C(t, x, u(t, x)) \Big) dtdx  
- \int_V \int_0^T \Big(X^{\circ}(t, x)\hat{A}^{\ast}p^{\circ}(t, x) - p^{\circ}(t, x)\hat{A}X^{\circ}(t, x) \Big)dtdx
~~~~~~\\
+ \int_V \Big(p^{\circ}(0, x)\xi(x) + \int_0^T p^{\circ}(t, x)\sigma(t, x)dB(t)\Big) dx \Big].~~~~~ \nonumber
\end{eqnarray}
On the other hand, since $\bar{X}(t, x)$ satisfies
\begin{eqnarray*}
d\bar{X}(t, x) = \big(\hat{A}\bar{X}(t, x) + C(t, x, u(t, x)\big)dt + \sigma(t, x)dB(t), 
\end{eqnarray*}
we get by integration of parts ([2])
\begin{eqnarray}
\int^T_0 \frac{dp^{\circ}(t, x)}{dt}\bar{X}(t, x)dt =
 p^{\circ}(T, x)\bar{X}(T, x) - p^{\circ}(0, x)\bar{X}(0, x)  ~~~~~~~~~~~~
~~~~~~~~~~~~~~~~~~~~~~~~~~~
 \nonumber \\ 
 - \int^T_0 p^{\circ}(t, x)
\big[ \hat{A}\bar{X}(t, x)dt + C(t, x, u(t, x))dt + \sigma(t, x)dB(t)\big]  ~~~~~~~~~~~~~~\\
~~~~~~~~~~~~~~~~~ =  - p^{\circ}(0, x)\xi(x) -  \int^T_0 p^{\circ}(t, x)
\big[ \hat{A}\bar{X}(t, x)dt + C(t, x, u(t, x)) dt+ \sigma(t, x)dB(t)\big], \nonumber 
\end{eqnarray}
where we used $p^{\circ}(T, x) = 0$ and $\bar{X}(0, x) = \xi(x)$.   
Since $\bar{X} - X^{\circ} = 0$ (see Eqs.(3) and (10)) and $p^{\circ} = 0$ on $\Gamma$, the surface of $V$, the first Green formula ([8, p.258]) implies
\begin{eqnarray*}
\int_V \big(\bar{X}(t, x) - X^{\circ}(t, x)\big) \hat{A}^{\ast}p^{\circ}(t, x) dx = 
\int_V \ p^{\circ}(t, x)\hat{A}\big(\bar{X}(t, x) - X^{\circ}(t, x)\big)dx.
\end{eqnarray*}
From this equality it follows that
\begin{eqnarray}
\int_V \big(X^{\circ}(t, x)\hat{A}^{\ast}p^{\circ}(t, x) - p^{\circ}(t, x)\hat{A}X^{\circ}(t, x)\big) dx
=
\int_V \big(\bar{X}(t, x)\hat{A}^{\ast}p^{\circ}(t, x) - p^{\circ}(t, x)\hat{A}\bar{X}(t, x)\big) dx.
\end{eqnarray}
\indent
Upon substituting Eqs.(15), (16) into (14), we see that the right-hand side of Eq.(14) (and (13)) is equal to zero. Hence
we can conclude that for each $u \in \mathcal {A}$ it follows that
\begin{eqnarray*}
   {\mathcal J}(u) \le L(u; X^{\circ}, p^{\circ}).
\end{eqnarray*}
Since $(X^{\circ}, p^{\circ}) \in {\mathcal B}$ is arbitrary, we have
\begin{eqnarray*}
   {\mathcal J}(u) \le \inf_{(X^{\circ},~p^{\circ}) \in {\mathcal B}}~L(u; X^{\circ}, p^{\circ}).
\end{eqnarray*}
The optimal value for the primal problem is $\sup_{u \in \mathcal A} {\mathcal J}(u)$ and it satisfies
\begin{eqnarray*}
   \sup_{u \in \mathcal A} {\mathcal J}(u) \le \sup_{u \in \mathcal A} \inf_{(X^{\circ},~p^{\circ}) 
   \in {\mathcal B}}~L(u; X^{\circ}, p^{\circ}).
\end{eqnarray*}
By a well-known inequality of game theory [3], we have
\begin{eqnarray}
   \sup_{u \in \mathcal A} {\mathcal J}(u) \le \sup_{u \in \mathcal A} ~\inf_{(X^{\circ},~p^{\circ}) \in {\mathcal B}} 
   L(u; X^{\circ}, p^{\circ}) \le
   \inf_{(X^{\circ},~p^{\circ}) \in {\mathcal B}}~ \sup_{u \in \mathcal A}~L(u; X^{\circ}, p^{\circ}).
\end{eqnarray}
In view of (6) we see that for each fixed $(X^{\circ}, p^{\circ}) \in \mathcal B$ the value
$\sup_{u \in \mathcal A}~L(u; X^{\circ}, p^{\circ})$ is identical to Eq.(7) of the dual problem, that is, 
${\mathcal L}(X^{\circ}, p^{\circ})$, which is to be minimized.  
Therefore, we obtain
\begin{eqnarray*}
   \sup_{u \in \mathcal A} {\mathcal J}(u) \le \inf_{(X,~p) \in {\mathcal B}}{\mathcal L}(X, p).
\end{eqnarray*}
This proves the weak duality theorem.\\
\\
\indent
The last inequality shows that each $(X,~p) \in {\mathcal B}$ provides an upper bound for the primal problem.  \\
\\
\\
{\large \bf 4. Strong Duality Theorem} \\  
\\
\indent 
In this sction we assume that a control process $\bar{u}$ satisfies a sort of the maximum principle of optimality such as 
in [4, Theorem 2.1]. Under the concavity of the function $F$ in Eq.(4), it entails the strong duality theorem.
More precisely, the corresponding solution 
$\bar{X} = X_{\bar{u}}$ of Eqs.(1)-(3) provides an optimal control for the dual problem and there is no duality gap; 
both extreme values (5) and (11) are exactly equal. \\
\\
{\bf Theorem 2.} 
{\it Suppose $\bar{X}$ is a solution of Eqs.{\rm (1)-(3)} for an admissible control $\bar{u} \in \mathcal A$, and that 
$\bar{p}$, together with this $\bar{X}$, is a solution of Eqs.{\rm (8)-(10)}. If $\bar{u} \in \mathcal A$ satisfies 
\begin{eqnarray}
   H(t, x, \bar{p}(t, x)) = G(t, x, \bar{u}(t, x)) + \bar{p}(t, x)C(t, x, \bar{u}(t, x)),~~{\it for~all}~(t, x) \in [0, T] \times V,
\end{eqnarray}
the function $H$ being defined by {\rm (6)}, then $\bar{u}$ is an optimal control of the primal problem and
$\bar{X}$ is that of the dual one. Moreover, there is no duality gap;
\begin{eqnarray*}
   \sup_{u \in \mathcal A} {\mathcal J}(u) = \inf_{(X,~p) \in {\mathcal B}}{\mathcal L}(X, p).
\end{eqnarray*}
}
\\
{\it Proof.} The proof is similar to that of Theorem 1. Let us put
\begin{eqnarray*}
\mathcal{J}(\bar{u}) = \mathbb{E} \Big[ \int_0^T \Big( \int_V \big(F(t, x, \bar{X}(t, x)) + 
G(t, x, \bar{u}(t,x)) \big)dx \Big)dt \Big].
\end{eqnarray*}
On the other hand, using (7) and (18), we have
\begin{eqnarray*}
\mathcal{L}(\bar{X}, \bar{p}) = \mathbb{E} \Big[ \int_V \int_0^T \Big(F(t, x, \bar{X}(t, x)) 
+ G(t, x, \bar{u}(t, x)) + \bar{p}(t, x)C(t, x, \bar{u}(t, x)) ~~~~~~~~~~~~\\
~~~~~~~- \bar{X}(t, x) \frac{\partial F(t, x, \bar{X}(t, x))}{\partial X}\Big) dtdx  
-  \int_V \int_0^T \Big(\bar{X}(t, x)\hat{A}^{\ast}\bar{p}(t, x) - \bar{p}(t, x)\hat{A}\bar{X}(t, x) \Big)dtdx  \\
+ \int_V \Big(\bar{p}(0, x)\xi(x) + \int_0^T \bar{p}(t, x)\sigma(t, x)dB(t)\Big) dx \Big]. 
\end{eqnarray*}
We evaluate the difference
\begin{eqnarray*}
~~~~~\mathcal{J}(\bar{u}) - \mathcal{L}(\bar{X}, \bar{p}) =  
\mathbb{E} \Big[ \int_V \int_0^T \Big (\frac{\partial F(t, x, \bar{X}(t, x))}{\partial X}\bar{X}(t, x)
- \bar{p}(t, x)C(t, x, \bar{u}(t, x)) \Big) dtdx  \nonumber \\
+  \int_V \int_0^T \Big(\bar{X}(t, x)\hat{A}^{\ast}\bar{p}(t, x) - \bar{p}(t, x)\hat{A}\bar{X}(t, x) \Big)dtdx  
~~~~~~~~~~~~~\\
- \int_V \Big(\bar{p}(0, x)\xi(x) + \int_0^T \bar{p}(t, x)\sigma(t, x)dB(t)\Big) dx \Big].  \nonumber
\end{eqnarray*}
Now it is easy to prove that the difference is equal to zero, using a similar calculation to the right-hand side of (13); 
   $\mathcal J(\bar{u}) = {\mathcal L}(\bar{X}, \bar{p})$.
Using Theorem 1 (weak duality), 
it follows that $\bar u$ is an optimal control 
for the primal problem and that $(\bar X, \bar p)$ is an optimal pair for the dual one.
This completes the proof.  \\
\\
\indent
Although our system is simpler than that of [4] and the approach is different from it, Eq.(18) turns out
a sufficient optimality condition for the primal problem. \\
\\
\\
{\large \bf 5. Partial Observation Control} \\  
\\
\indent 
In partially observable systems as in [1], it is necessary to consider controls that do not depend on the space 
variable $x$. We denote the subset of such controls by $\mathcal A_1$; 
\begin{eqnarray*}
\mathcal A_1 = \{ u(t) = u(t, \omega) ~|~ u(t) \in \mathcal A \}.
\end{eqnarray*}
The primal problem is to maximize the functional 
\begin{eqnarray*}
{\mathcal J}(u) = \mathbb{E} \Big[ \int_0^T \Big( \int_V \big(F(t, x, X(t, x)) + G(t, x, u(t)) \big)dx \Big)dt \Big] 
\end{eqnarray*}
over $u \in \mathcal A_1$ together with $X(t, x)$ satisfying
\begin{eqnarray*}
dX(t, x) & = & \big(\hat{A}X(t, x) + C(t, x, u(t))\big)dt + \sigma(t, x) dB(t), \\
X(0, x) & = & \xi(x)~ ~{\rm for}~x \in \bar{V},~~
X(t, x) =  \eta(t, x)~~{\rm for}~(t, x) \in (0, T) \times \Gamma. 
\end{eqnarray*}
\indent
The dual system is governed by Eqs.(8)-(10) as before. In order to formulate the dual problem, let us put
\begin{eqnarray}
\mathcal{H}(t, p(t, \cdot)) = \sup_{u \in \mathcal U} \int_V \big(G(t, x, u) + p(t, x)\,C(t, x, u) \big) dx,
\end{eqnarray}
for functions $p(t, x)$ that are solutions of Eqs.(8)-(10). The dual problem is to minimize the functional
\begin{eqnarray*}
\mathcal{L}_1(X, p) = \mathbb{E} \Big[ \int_0^T \mathcal{H}(t, p(t, \cdot))dt \Big]  
+ \mathbb{E} \Big[ \int_V \int_0^T \Big(F(t, x, X(t, x)) - 
X(t, x) \frac{\partial F(t, x, X(t, x))}{\partial X}\Big) dtdx  \nonumber\\
- \int_V \int_0^T \Big(X(t, x)\hat{A}^{\ast}p(t, x) - p(t, x)\hat{A}X(t, x) \Big)dtdx  \\
+ \int_V \Big(p(0, x)\xi(x) + \int_0^T p(t, x)\sigma(t, x)dB(t)\Big) dx \Big],   
\end{eqnarray*}
over all $p(t, x)$ and $X(t, x)$ satisfying Eqs.(8)-(10). 
Note that this type of dual problem takes a more similar form to the one dealt with in [7]. \\
\indent
We prove two duality theorems. To do this, let us take an arbitrarily chosen control $u \in \mathcal A_1$, and 
introduce the corresponding functional $L(u; X, p)$ similar to Eq.(12), while $X, p$
satisfy Eqs.(8)-(10), {\it i.e.,} $(X, p) \in \mathcal B$. Then we can derive the inequality analogous to (17);
\begin{eqnarray*}
   \sup_{u \in \mathcal A_1} {\mathcal J}(u)  \le
   \inf_{(X,~p) \in {\mathcal B}}~ \sup_{u \in \mathcal A_1}~L(u; X, p),
\end{eqnarray*}
for all $(X, p) \in \mathcal B$. 
Among the terms of $L(u; X, p)$, those relevant to $u(t)$ are $G(t, x, u(t))$ and $p(t, x)C(t, x, u(t))$.
Hence we divide $\sup_{u \in \mathcal A_1} L(u; X, p)$ into two parts: one is
\begin{eqnarray}
\sup_{u \in \mathcal A_1} \mathbb{E} \Big[ \int_0^T \big(\int_V (G(t, x, u(t)) + p(t, x)C(t, x, u(t)))dx\big)dt \Big];
\end{eqnarray}
and the other is
\begin{eqnarray}
\mathbb{E} \Big[ \int_V \int_0^T \Big(F(t, x, X(t, x)) - X(t, x) \frac{\partial F(t, x, X(t, x))}{\partial X}\Big) dtdx ~~~~~~
~~~~~~~~~~~~~~~~~~~~~~~~~~~~~~ 
\end{eqnarray}
\begin{eqnarray*}
- \int_V \int_0^T \Big(X(t, x)\hat{A}^{\ast}p(t, x) - p(t, x)\hat{A}X(t, x) \Big)dtdx  
+ \int_V \Big(p(0, x)\xi(x) + \int_0^T p(t, x)\sigma(t, x)dB(t)\Big) dx \Big].  
\end{eqnarray*}
Using a measurable selection theorem, Fubini's theorem and Eq.(19), we see that the expectation (20) can be written as
\begin{eqnarray*}
\mathbb{E} \Big[ \int_0^T \mathcal{H}(t, p(t, \cdot))dt \Big].
\end{eqnarray*}
This together with (21) yields the functional $\mathcal{L}_1(X, p)$ for which the weak duality theorem holds;  
\begin{eqnarray*}
\sup_{u \in \mathcal A_1} {\mathcal J}(u)  \le
   \inf_{(X,~p) \in {\mathcal B}}~\mathcal{L}_1(X, p).
\end{eqnarray*}
\indent
Next suppose that $\bar{X}$ is a solution of Eqs.(1)-(3) for an admissible control $\bar{u} \in \mathcal A_1$, and that
$(\bar{X}, \bar{p}) \in \mathcal B$ satisfies (averaged maximum condition in [4])
\begin{eqnarray*}
\mathcal{H}(t, \bar{p}(t, \cdot)) = \int_V \big(G(t, x, \bar{u}(t)) + \bar{p}(t, x)\,C(t, x, \bar{u}(t)) \big) dx
~~~{\rm ~for ~all}~ t \in [0, T].
\end{eqnarray*}
Then we obtain the equality ${\mathcal J}(\bar{u}) = \mathcal{L}_1(\bar{X}, \bar{p})$ and hence 
the strong duality theorem as in Section 4, implying no duality gap
\begin{eqnarray*}
\sup_{u \in \mathcal A_1} {\mathcal J}(u) = \inf_{(X,~p) \in {\mathcal B}}~\mathcal{L}_1(X, p).
\end{eqnarray*}
Moreover, from the weak duality theorem it follows that $\bar{u}$ provides an optimal control for the primal problem, and so does the pair $(\bar{X}, \bar{p})$ 
for the dual problem.\\ 
\\
\\
{\bf References}
\begin{itemize}
\item[{[1]}]   A. Bensoussan, {\it Stochastic Control of Partially Observable Systems}, Cambridge Univ. Press, 1992.
\item[{[2]}]   K. L. Chung and R. J. Williams, {\it Introduction to Stochastic Integration}, Second Edition, 
Birkh\"{a}user, 1990.
\item[{[3]}]   S. Karlin, {\it Mathematical Methods and Theory in Games, Programming and Economics}, Vol. I,
Addison-Wesley, 1959.
\item[{[4]}]  B. {\O}ksendal, Optimal control of stochastic partial differential equations,
{\it Stochastic Analysis and Applications}, {\bf 23}, No. 1, 165--179, 2005.
\item[{[5]}]  S. Tanimoto, A duality theorem for max-min control problems,
{\it IEEE Transactions on Automatic Control}, {\bf AC-27}, No. 5, 1129--1131, 1982.
\item[{[6]}]  S. Tanimoto, Duality in the optimal control of non-well-posed distributed systems,
{\it Journal of Mathematical Analysis and Applications}, {\bf 171}, 277--287, 1992.
\item[{[7]}]  S. Tanimoto, Duality and lower bounds in optimal stochastic control,
{\it International Journal of Systems Science}, {\bf 25}, 1365--1372, 1994.
\item[{[8]}]  J. Wloka, {\it Partial Differential Equations}, Cambridge Univ. Press, 1987.
\end{itemize}
\end{document}